\documentclass[11pt,twoside,leqno,draft]{article}

\usepackage{amsmath,amssymb,amsthm,bm}
\numberwithin{equation}{section}

\theoremstyle{plain}
\newtheorem{thm}{Theorem}[section]
\newtheorem{prop}[thm]{Proposition}
\newtheorem{lem}[thm]{Lemma}
\newtheorem{cor}[thm]{Corollary}
\theoremstyle{definition}

\newtheorem{example}{Example}[section]
\theoremstyle{remark}
\newtheorem{remark}{Remark}[section]

\let \dummy=\epsilon
\let \epsilon=\varepsilon
\let \varepsilon=\dummy
\let \dummy=\phi
\let \phi=\varphi
\let \varphi=\dummy
\let \dummy=\relax

\newcommand{\el}{l}
\newcommand{\bC}{\mathbb{C}}
\newcommand{\bR}{\mathbb{R}}
\newcommand{\bZ}{\mathbb{Z}}
\newcommand{\bZp}{\bZ_{>0}}
\newcommand{\bZnn}{\bZ_{\ge0}}
\newcommand{\abs}[1]{\left|#1\right|}
\newcommand{\deq}{:=}
\newcommand{\prims}[1][p]{#1\text{\upshape\,:\,prime}}
\newcommand{\qbinom}[3][q]
{\genfrac{[}{]}{0pt}{}{#2}{#3}_{#1}}
\newcommand{\pbinom}[3][p^{-s}]
{\genfrac{[}{]}{0pt}{}{#2}{#3}_{#1}}
\newcommand{\floor}[1]{\lfloor#1\rfloor}
\newcommand{\length}[1]{\ell(#1)}
\newcommand{\term}[1]{{\itshape #1}}
\newcommand{\textem}[1]{{\bfseries #1}}
\newcommand{\by}[2]{#1, #2}  
\newcommand{\volume}[1]{{\bfseries #1}}
\newcommand{\journal}[1]{{\itshape #1}}
\bmdefine{\blambda}{\lambda}

\let \Re=\relax
\let \Im=\relax
\DeclareMathOperator{\Re}{\mathrm{Re}}
\DeclareMathOperator{\Im}{\mathrm{Im}}
\DeclareMathOperator{\ord}{ord}
\DeclareMathOperator{\Mult}{Mult}
\DeclareMathOperator*{\Res}{Res}

\title{{\bfseries Multiple finite Riemann zeta functions}}
\author{Kazufumi KIMOTO, Nobushige KUROKAWA,\\
Sho MATSUMOTO and Masato WAKAYAMA}
\date{\today}
\begin{document}
\setlength{\baselineskip}{18pt}
\maketitle

\begin{abstract}
Observing a multiple version of the divisor function
we introduce a new zeta function
which we call a \term{multiple finite Riemann zeta function}.
We utilize some $q$-series identity for proving the zeta function has an 
Euler product and then, describe the location of zeros. 
We study further 
multi-variable and multi-parameter versions  of the 
multiple finite Riemann 
zeta functions and their {\it infinite} counterparts in connection with 
symmetric polynomials and 
some arithmetic quantities called \term{powerful numbers}. 

\par\noindent
\textem{2000 Mathematics Subject Classification} : 11M36
\par\noindent
\textem{Key Words} : Riemann zeta function, divisor functions, 
Euler product, functional equations, $q$-binomial coefficients, 
complete symmetric polynomial, powerful numbers, Eisenstein series.
\end{abstract}

\section{Introduction}

The divisor function $\sigma_k(N)\deq\sum_{d|N}d^k$
is a basic multiplicative function
and plays an important role 
from the beginning of the modern arithmetic study,
and in particular,
$\sigma_k(N)$ appears in the Fourier coefficients
of the (holomorphic) Eisenstein series $E_{k+1}(\tau)$.
Usually perhaps,
it is not so common to regard $\sigma_k(N)$
as a sort of zeta function.
However, in the present paper
we treat the divisor function $\sigma_k(N)$
as a function of a complex variable $k=-s$. 
Still, it is clear that when $N\to\infty$ through factorials (or $0$) we have 
$\sigma_{-s}(N) \to \zeta(s)$, the Riemann zeta function, 
and there are at least two interpretations of 
$$
Z_N^1(s):=\sum_{n|N} n^{-s} =\sigma_{-s}(N)
$$
as a zeta function in number theory:

$\bullet$ Fourier coefficients of real analytic Eisenstein series,

$\bullet$ Igusa zeta functions.

Concerning the first, we refer to Bump et al \cite{BCKV}
where the so-called ``local Riemann hypothesis" is studied. In case of
the real analytic Eisenstein series $E(s,\tau)$ for the modular group
$SL(2,\mathbb Z)$, the $N$-th Fourier coefficient is essentially given 
by
$$
c_N(s,\tau):
=Z_N^1(2s-1)K_{s-\frac12}(2\pi N \Im(\tau))e^{2\pi i
\Re(\tau)}
$$
and hence it satisfies the local Riemann hypothesis:
if $c_N(s,\tau)=0$ then $\Re(s)=\frac12$. 
On the other hand, for the second, the interpretation is coming from the
(global) 
Igusa zeta function $\zeta^{\textrm{Igusa}}(s, R)$ of 
a ring $R$ defined as 
$$
\zeta^{\textrm{Igusa}}(s, R) \deq \sum_{m=1}^\infty
{}^\#\textrm{Hom}_{\textrm{ring}} (R, \mathbb Z/(m)) m^{-s}.
$$
Then, in fact, it is easy to see that 
$$
Z_N^1(s)=\zeta^{\textrm{Igusa}}(s, \mathbb Z/(N)).
$$
 
Thus, the purpose of the present 
paper is initially to study the function defined by the series 
\begin{equation}\label{MMzeta}
Z^{m}_N(s)\deq\sum_{n_1|n_2|\cdots|n_m|N}
(n_1n_2\cdots n_m)^{-s}.
\end{equation}
We call $Z^{m}_N(s)$
the {multiple finite Riemann zeta function}
of type $N$.
We then study several basic properties of $Z^{m}_N(s)$
such as an Euler product, a functional equation and
an analogue of the Riemann hypothesis in an elementary way
by the help of some $q$-series identity.
We also study the limit case $Z^{m}_\infty(s)\deq \lim_{N\to\infty}Z^{m}_N(s) (=Z^m_0(s))
=\sum_{n_1|n_2|\cdots|n_m}
(n_1n_2\cdots n_m)^{-s}$.
Here and throughout the paper we consider the limit $N \to \infty$ as the one through factorials. 

Moreover,
we generalize this zeta function in two directions.
The first one is to increase the number of variables.
We prove that the Euler product 
of a multi-variable version of the zeta function is expressed 
in terms of the complete symmetric polynomials
with a remarkable specialization of variables 
(Theorem \ref{CompleteSymmetric}).
The second one is further to add parameters indexed by a set of 
positive integers.
For general parameters 
it seems difficult to calculate an explicit expression
of the Euler product using symmetric functions with some meaningful  
specialization of variables.
However, when we restrict ourselves to the one variable case, 
under a special but non-trivial specialization of parameters, 
we show that the corresponding multiple zeta functions
are written as a product of Riemann's zeta functions $\zeta(cs)$ with several 
constants $c$'s determined by the given parameters 
and the Dirichlet series associated with 
generalized powerful numbers (see Section 4 for the definition).
Moreover, we determine the condition whether  
the Dirichlet series associated with such generalized powerful numbers
can be extended as a meromorphic function
to the entire complex plane $\bC$ or not (see Theorem 4.8 and its corollary).
As a consequence, the most of 
such Dirichlet series are shown to have the imaginary axis as 
a natural boundary. The result is a generalization of the one in 
\cite{IvicShiu}.

In the final position of the paper,
we make a remark on the relation between
the multiple finite Riemann zeta functions and
the number of isomorphism classes of abelian groups.
In particular, using the Tauberian theorem, for a given positive integer $m$,  
we prove that the asymptotic average of the isomorphism classes
of abelian groups of order $n$ 
which are given by the direct sum of $p$-groups $A_p$
such that $p^m A_p=0$,
equals $\zeta(2) \zeta(3) \dotsb \zeta(m)$ (when $n$ tends to $\infty$). 
We make also a small discussion on an analogous notion of (holomorphic) 
multiple Eisenstein series
which are defined via these multiple finite Riemann zeta functions
as coefficients of its  Fourier series 
(or related to the so-called Lambert type series).

Throughout the paper, we denote 
the sets of all integers, positive integers,
non-negative integers, real numbers and complex numbers
by $\bZ$, $\bZp$, $\bZnn$, $\bR$ and $\bC$
respectively.

\section{Multiple finite Riemann zeta function}

In this section we prove the fundamental properties of
the multiple finite Riemann zeta function $Z^{m}_N(s)$
defined by \eqref{MMzeta} and make a discussion on 
some related Dirichlet series.
Namely we first show the following theorem.

\begin{thm}\label{Theorem1}
Let $N$ be a positive integer.
\begin{enumerate}
\item Euler product {\upshape:}
\begin{equation}\label{EulerProduct}
Z^{m}_N(s)=\prod_{\prims}\prod_{k=1}^m
\frac{1-p^{-s(\ord_pN+k)}}{1-p^{-sk}},
\end{equation}
where $\ord_pN$ denotes the order of 
$p$-factor in the prime decomposition of $N$. 
\item Functional equation {\upshape:}
\begin{equation}\label{FEq}
Z^{m}_N(-s)=N^{ms}Z^{m}_N(s).
\end{equation}
\item Analogue of the Riemann hypothesis {\upshape:}
All zeros of $Z^{m}_N(s)$ lie on the imaginary axis $\Re s=0$.
More precisely, the zeros of $Z^{m}_N(s)$ are
of the form $s=\frac{2\pi in}{(\ord_pN+k)\log p}$
for $k=1,\ldots,m,\, p|N$ and $n\in\bZ\setminus\{0\}$.
Consequently, the order $\Mult^m(n,p,k)$ of the zero
at $s=\frac{2\pi in}{(\ord_pN+k)\log p}$
is given by $\Mult^m(n,p,k)
=\#\{(\el,j)\;;\; 1\leq \el \leq m,\, j\in \bZ\setminus\{0\},\;
(\ord_pN+k)j=(\ord_pN+\el)n \}$.
\item Special value {\upshape:}
When $n$ is a positive integer one has $Z_N^m(-n)\in\bZ$.
\end{enumerate}
\end{thm}

For the proof of the theorem the following lemma is crucial.

\begin{lem}\label{SimpleFormulas}
Let $m$ be a positive integer. Then 
\begin{enumerate}
\item For any integers $\el\geq0$ it holds that
\begin{equation}\label{SkewBinom}
\sum_{d=0}^\el \qbinom{m-1+d}{m-1} q^d = \qbinom{m+ \el}{m},
\end{equation}
where $\qbinom{n}{k}$ is the $q$-binomial coefficient defined by
$\qbinom{n}{k} = \prod_{j=1}^k (1-q^{n+1-j})/(1-q^j)$.
\item 
For $\abs{x}<1, \, \abs{q}<1$ it holds that
\begin{equation}\label{SkewBinomInfty}
\sum_{d=0}^\infty \qbinom{m+d}{m} x^d
=\prod_{k=0}^m \frac{1}{1-q^k x}.
\end{equation}
\end{enumerate}
\end{lem}

\begin{proof}
We prove the formula \eqref{SkewBinom} by induction on $\el$.
When $\el=0$ the formula \eqref{SkewBinom} clearly holds. 
Suppose that it holds for $\el$.
Then we see that
\begin{align*}
& \sum_{d=0}^{\el+1} \qbinom{m-1+d}{m-1} q^d =
\qbinom{m+ \el}{m} + \qbinom{m+ \el}{m-1} q^{\el+1} \\
=& \prod_{k=1}^m (1-q^k)^{-1} \prod_{k=1}^{m-1} (1-q^{\el+1+k})
\left\{ (1-q^{\el+1} ) + (1-q^m) q^{\el+1} \right\} \\
=& \prod_{k=1}^m (1-q^k)^{-1} \prod_{k=1}^m (1-q^{\el+1+k})
= \qbinom{m+\el+1}{m},
\end{align*}
whence \eqref{SkewBinom} is also true for $\el+1$.

The second formula \eqref{SkewBinomInfty} can be proved in
the same manner by induction, but on $m$.
(It is also obtained from the so-called $q$-binomial theorem. 
See, e.g. \cite{AAR:SF}).
This completes the proof of the lemma.
\end{proof}

\begin{proof}[Proof of Theorem \ref{Theorem1}]
The functional equation \eqref{FEq} is easily seen
from the definition. Actually, we have 
\begin{align*}
Z^{m}_N(-s)=\sum_{n_1|n_2|\dotsb|n_m|N}
(n_1n_2\cdots n_m)^{s}
=N^{ms}\sum_{n_1|n_2|\dotsb|n_m|N}
\left(\frac{n_1}N\frac{n_2}N \dotsb \frac{n_m}N\right)^{s}
=N^{ms}Z^{m}_N(s).
\end{align*}

Next we prove the Euler product \eqref{EulerProduct}.
We first show that $Z_N^m(s)$ is multiplicative
with respect to $N$.
In fact, if we suppose $N$ and $M$ are co-prime,
then we observe
\begin{align*}
Z_{NM}^m(s)
&=\sum_{n_1|n_2|\dotsb|n_m|NM}(n_1n_2\dotsb n_m)^{-s}\\
&=\sum_{c_1|c_2|\dotsb|c_m|N}(c_1 c_2\dotsb c_m)^{-s}
\sum_{d_1|d_2|\cdots|d_m|M}(d_1d_2\cdots d_m)^{-s}
=Z_N^m(s)Z_M^m(s),
\end{align*}
because every divisor $n$ of $NM$ is uniquely written
as $n= c d$  where $c|N$ and $d|M$.
By means of this fact,
in order to get the Euler product expression \eqref{EulerProduct}
of $Z_N^m(s)$ it suffices to calculate
the case where $N$ is a power of prime $p$.
In this case, one proves the expression
\begin{align}\label{pEuler}
Z_{p^\el}^m(s)
=\prod_{k=1}^m\frac{1-p^{-(\el+k)s}}{1-p^{-sk}}
\end{align}
by induction on $m$ as follows:
It is clear that \eqref{pEuler} holds for $m=1$.
Now, suppose \eqref{pEuler} is true for $m-1$.
Then we see that
\begin{align*}
Z_{p^\el}^m(s)
&=\sum_{0\leq j_1\leq\dots\leq j_m\leq\el}
p^{-(j_1+\cdots+j_m)s}\\
&=\sum_{j_m=0}^\el Z_{p^{j_m}}^{m-1}(s)p^{-j_ms}
=\sum_{d=0}^\el\prod_{k=1}^{m-1}
\frac{1-p^{-(d+k)s}}{1-p^{-sk}}p^{-ds}
=\sum_{d=0}^\el \pbinom{m-1+d}{d} p^{-ds}.
\end{align*}
Therefore the assertion follows immediately
from the formula \eqref{SkewBinom} in the lemma above.
This proves \eqref{pEuler},
whence the Euler product for $Z_N^m(s)$ follows.

Using the Euler product \eqref{EulerProduct},
we observe that the meromorphic function $Z_N^m(s)$ may have zeros
at each $s=\frac{2\pi in}{(\ord_pN+k)\log p}$
for $k=1,\ldots,m,\; p|N$ and $n\in\bZ$.
Note however, since
\begin{equation}\label{0-value}
Z_N^m(0)=\prod_{\prims}\binom{\ord_pN+m}{m}\neq0,
\end{equation}
$s=0$ is not a zero of $Z_N^m(s)$.
Suppose now
\begin{equation*}
\frac{2\pi in}{(\ord_pN+k)\log p}
=\frac{2\pi im}{(\ord_qN+\el)\log q}
\end{equation*}
holds for some
$1\leq k, \el\leq m,\; p,q|N$ and $n,j\in \bZ\setminus\{0\}$.
Then it becomes $p^{(\ord_pN+k)j}=q^{(\ord_qN+\el)n}$.
This immediately shows that $p=q$ and
$(\ord_pN+k)j=(\ord_pN+\el)n$.
Hence the order of zero at $s=\frac{2\pi in}{(\ord_pN+k)\log p}$
is actually given by $\Mult^m(n,p,k)$.
Obviously, one has $\Mult^1(n,p,1)=1$.

The last claim about the special value of $Z_N^m(s)$ 
is clear from the definition.
This completes the proof of the theorem.
\end{proof}

\begin{remark}
Note that
$Z^m_{p^\el}(s) = Z^{\el}_{p^m}(s)$ from \eqref{pEuler}. \qed
\end{remark}

Letting $N\to \infty$ (or $N\to 0$) in Theorem \ref{Theorem1},
we have the following 

\begin{cor}\label{MMzetaInfty}
Define a function $Z^m_\infty(s)$ by
\begin{equation*}
Z^m_\infty(s)\deq\sum_{\begin{subarray}{c}
n_1|n_2|\cdots|n_m \\
n_i\in\bZp
\end{subarray}}
(n_1n_2\cdots n_m)^{-s}  \qquad (\Re s>1).
\end{equation*}
Then we have $Z^m_\infty(s)=\prod_{k=1}^m \zeta(ks)$.
Here $\zeta(s)$ is the Riemann zeta function.
\qed
\end{cor}

We next consider a Dirichlet series defined
via the multiple finite Riemann zeta functions.

\begin{cor}
Retain the notation in the corollary above.
Define a Dirichlet series $\zeta^m(s)$ by
\begin{equation*}
\zeta^m(s)\deq\sum_{n=1}^\infty Z_n^m(s)n^{-s}.
\end{equation*}
Then we have
\begin{equation}\label{m-pleZeta}
\zeta^m(s)=\prod_{k=1}^{m+1}\zeta(sk)=Z^{m+1}_\infty(s).
\end{equation}
\end{cor}

\begin{proof}
Since $Z_n^m(s)$ is multiplicative relative to $n$,
the Dirichlet series $\zeta^m(s)$ has an Euler product of the form
\begin{equation*}
\zeta^m(s)=\prod_{\prims}\left(\sum_{\el=0}^\infty
Z_{p^\el}^m(s)p^{-s\el}\right).
\end{equation*}
Hence the first equality follows from Lemma \ref{SimpleFormulas}.
In fact,
\begin{align*}
\sum_{\el=0}^\infty Z_{p^\el}^m(s)p^{-s\el}
=\sum_{\el=0}^\infty \prod_{k=1}^m
\frac{1-p^{-s(\el+k)}}{1-p^{-sk}} p^{-s\el}
=\prod_{k=1}^{m+1} \frac1{1-p^{-sk}}
\end{align*}
when $\Re s>0$. 
The second equality is obvious from Corollary \ref{MMzetaInfty}.
This proves the corollary.
\end{proof}

\begin{remark}
Let $g(n)$ be the number of the isomorphism classes
of the abelian groups of order $n$.
Then we have (see, e.g. \cite[p.16]{Z})
\begin{equation*}
\lim_{m\to\infty}\zeta^m(s)=\prod_{k=1}^\infty \zeta(ks)
=\sum_{n=0}^\infty g(n)n^{-s}
\qquad (\Re s>1).
\end{equation*}
Observing the analytic property of $Z^m_N(s)$, 
we have much precise information about the isomorphism classes 
and their asymptotic averages. 
For details, see Section \ref{ConcludingRemarks}.
\qed
\end{remark}

We now slightly generalize the discussion above. 
Put
\begin{equation*}
\zeta^m(s,t)\deq\sum_{n=1}^\infty Z_n^m(s) n^{-t}.
\end{equation*}
It is clear that $\zeta^m(s)=\zeta^m(s,s)$.
Then, by the same reasoning,
$\zeta^m(s,t)$ has the Euler product
\begin{equation}\label{m-pleZeta(s,t)}
\zeta^m(s,t)
=\prod_{\prims}
\left(\sum_{\el=0}^\infty Z_{p^{\el}}^m(s) p^{-\el t}\right)
=\prod_{k=0}^m\zeta(sk+t).
\end{equation}
Actually, it follows from the formula \eqref{SkewBinomInfty} that
\begin{equation*}
\sum_{\el=0}^\infty
\left(\prod_{k=1}^m
\frac{1-p^{-s(k+\el)}}{1-p^{-sk}}p^{-t \el}\right)
= \prod_{k=0}^m \frac1{1-p^{-sk-t}}
\end{equation*}
for $\Re s >0$ and $\Re t >0$.
This is a generalization of the formula \eqref{m-pleZeta}.
One may consider $\zeta^m(s,t)$ as a zeta function for two variables.
We give here a few examples:

\begin{example}
The following is quite well-known.
For $\Re t> k+1$, we have
\begin{align*}
\zeta^1(-k,t)=\sum_{n=1}^\infty Z_n^1(-k)n^{-t}
=\sum_{n=1}^\infty \sigma_{k}(n)n^{-t}=\zeta(t)\zeta(t-k).
\end{align*}
\qed
\end{example}

\begin{example}
From \eqref{0-value},
for $\Re t>1$, it follows that 
\begin{equation*}
\zeta^m(0,t)  
= \sum_{n=1}^\infty Z_n^m(0) n^{-t}
= \sum_{n=1}^{\infty} \prod_{\prims} \binom{m+\ord_p n}{m} n^{-t}
= \zeta(t)^{m+1}.
\end{equation*}
\qed
\end{example}

\begin{example}
Let $\el$ be a positive integer.
For $\Re t>1-\el$, we have
\begin{equation*}
\zeta^m(\el, t+\el) = \prod_{k=1}^{m+1} \zeta(t+ \el k).
\end{equation*}
As $m \to \infty$,
$\lim_{m \to \infty} \zeta^m(\el,t+\el)
=\prod_{k=1}^{\infty} \zeta(t+\el k)$
is the higher Riemann zeta function studied in \cite{KMW}
(see also \cite{KW1}). It should be noted that the higher Riemann zeta function 
$\zeta_{\el \infty}(s):=\prod_{k=1}^{\infty} \zeta(s+\el k)$ possesses a functional equation. 
\qed
\end{example}

\section{Multivariable version}

We generalize the definition \eqref{MMzeta}
of the multiple finite Riemann zeta functions $Z^{m}_N (s)$. 
For $\gamma=(\gamma_1, \dots, \gamma_m) \in \bZp^m$ and $N \in \bZp$, 
we define
\begin{equation}
Z^{\gamma}_N (t_1, \dots, t_m)
\deq\sum_{ n_m^{\gamma_m}| \cdots | n_1^{\gamma_1}|N}
n_1^{- \gamma_1 t_1} \cdots n_m^{-\gamma_m t_m}.
\end{equation}
This is multiplicative with respect to $N$.
We also notice that
$Z^m_N (s) = Z^{\gamma}_N(s, \dots, s)$
when $\gamma = (1^m) =(\overbrace{1, \dots, 1}^m)$.
We can prove the following lemma
in a similar way to Theorem \ref{Theorem1}. 

\begin{lem} 
For $\el \geq 0$,
define a function $G^{\gamma}_{\el} (q_1, q_2, \dots, q_m)$ by
\begin{equation}\label{generating}
G^{\gamma}_{\el}(q_1, \dots, q_m) 
= \sum_{\begin{subarray}{c}
\el \geq \lambda_1 \geq \lambda_2 \geq \dots \geq \lambda_m \geq 0 \\
\gamma_j | \lambda_j
\end{subarray}}
q_1^{\lambda_1} \cdots q_m^{\lambda_m},
\end{equation}
where the sum is taken over all partitions
$\lambda=(\lambda_1, \lambda_2, \dots)$ of length $\leq m$
such that
$\lambda_1 \leq \el$ and $\gamma_j| \lambda_j$ for $ 1 \leq j \leq m$.
Then the Euler product of $Z^{\gamma}_N (t_1, \dots, t_m)$ is given as  
\begin{equation*}
Z^{\gamma}_N(t_1, \dots, t_m)
= \prod_{\prims}
G^{\gamma}_{\ord_p N} (p^{-t_1}, \dots, p^{-t_m}).
\end{equation*}
\qed
\end{lem}

\begin{remark} Let $f(n_1,\dots, n_m)$ be a function defined on $\in \bZp^m$ 
which is multiplicative with respect to each variable $n_j$. 
Then, in general, we can define a multiple zeta function by 
$$
Z^{\gamma}_N(f)
= \sum_{n_m^{\gamma_m} | \cdots | n_1^{\gamma_1} |N}
f(n_1, \dots, n_m). 
$$
Actually, one can show that 
$Z^{\gamma}_N(f)$ is multiplicative with respect to $N$, whence 
has the Euler product as 
$Z^{\gamma}_N(f)
= \prod_{\prims} Z^{\gamma}_{p^{\ord_p N}} (f)$.
\qed
\end{remark}

We first look at the simplest case where $\gamma=(1^m)$.
We abbreviate  respectively 
$Z^{(1^m)}_N(t_1, \dots, t_m)$
and $G^{(1^m)}_{\el}(q_1, \dots, q_m)$
to $Z^m_N(t_1, \dots, t_m)$
and $G^m_{\el}(q_1, \dots, q_m)$.

\begin{thm} \label{CompleteSymmetric}
Let $h_j(x_1, \dots, x_m)$ be
the $j$-th complete symmetric polynomial defined by 
\begin{equation*}
h_j(x_1, \dots, x_m)
= \sum_{\begin{subarray}{c}
i_1+ \cdots +i_m=j \\
i_k \in \bZnn
\end{subarray}}
x_1^{i_1} \cdots x_m^{i_m}.
\end{equation*}
Then we have
\begin{equation}\label{complete}
G^m_{\el} (q_1, \dots, q_m)
= \sum_{j=0}^{\el} h_j(q_1, q_1 q_2, \dots, q_1 q_2 \cdots q_m).
\end{equation}
In particular, we have
\begin{equation}\label{etatypeproduct}
G^m_{\infty} (q_1, \dots, q_m)
\deq\lim_{\el \to \infty} G^m_{\el}(q_1, \dots, q_m)
= \prod_{k=1}^{m} \frac1{1-q_1 q_2 \cdots q_k}.
\end{equation}
\end{thm}

\begin{proof}
Since
\begin{align*}
h_j(q_1, q_1 q_2, \dots, q_1 q_2 \cdots q_m)
&= \sum_{i_1+i_2+ \cdots + i_m =j}
q_1^{i_1+i_2+ \cdots +i_m}  q_2^{i_2+ \cdots + i_m} \cdots q_m^{i_m} \\
&= \sum_{j = \lambda_1 \geq \lambda_2 \geq \dots \geq \lambda_m \geq 0}
q_1^{\lambda_1} q_2^{\lambda_2} \cdots q_m^{\lambda_m},
\end{align*}
the first formula \eqref{complete} is clear 
from the definition \eqref{generating}.
The second formula \eqref{etatypeproduct} follows 
from \eqref{complete} together with the fact that 
the generating function of complete symmetric polynomials is given by 
$\sum_{j=0}^{\infty} h_j(x_1, \dots, x_m) z^j
= \prod_{k=1}^m (1-x_k z)^{-1}$
(see \cite{Mac}).
\end{proof}

As a corollary of the theorem, we obtain the Euler product of 
$Z^m_N (t_1, \dots, t_m)$.
\begin{cor}
Let $m$ and $N$ be positive integers.
Then we have
\begin{equation*}
\begin{split}
Z^m_N (t_1, \dots, t_m)
&= \sum_{n_m | \cdots | n_1|N} n_1^{-t_1} \cdots n_m^{-t_m}\\
&= \prod_{\prims} \left(
\sum_{j=0}^{\ord_p N}
h_j(p^{-t_1}, p^{-t_1-t_2}, \dots, p^{-t_1-t_2- \cdots -t_m})
\right).
\end{split}
\end{equation*}
Further, for $\Re t_j >1\ (1\leq j \leq m)$,
it follows that
\begin{equation*}
Z^m_{\infty}(t_1, \dots, t_m)
\deq\sum_{n_m |\cdots |n_1} n_1^{-t_1} \cdots n_m^{-t_m}
= \prod_{k=1}^m \zeta( t_1 + t_2 + \cdots + t_k).
\end{equation*}
\qed
\end{cor}

The following lemma gives a recurrence equation 
among $G^\gamma_\el(q_1,\dots,q_m)$'s.

\begin{lem}\label{LemmaA}
For $\gamma=(\gamma_1 ,\dots, \gamma_m) \in \bZp^m$, we have
\begin{equation*}
G^{\gamma}_{\el} (q_1, \dots, q_m) 
= \sum_{n=0}^{\floor{\el/\gamma_1}} q_1^{\gamma_1 n} 
G^{(\gamma_2, \dots, \gamma_m)}_{\gamma_1 n} (q_2, \dots, q_m). 
\end{equation*}
Here $\floor{x}$ is the largest integer not exceeding $x$.  
\end{lem}

\begin{proof}
Observe that
\begin{equation*}
G^{\gamma}_{\el}(q_1, \dots, q_m)
= \sum_{ 0 \leq n \leq \el/\gamma_1} q^{\gamma_1 n}
\sum_{\begin{subarray}{c}
\gamma_1 n \geq \lambda_2 \geq \dots \geq \lambda_m \geq 0, \\
\gamma_j |\lambda_j \ (2 \leq j \leq m)
\end{subarray}}
q_2^{\lambda_2} \cdots q_m^{\lambda_m} 
= \sum_{n=0}^{ \floor{\el/\gamma_1}} q_1^{\gamma_1 n} 
G^{(\gamma_2, \dots, \gamma_m)}_{\gamma_1 n} (q_2, \dots, q_m). 
\end{equation*}
\end{proof}

The following lemma shows that
it is enough to study the case where
$\gamma_1, \dots, \gamma_m$ are relatively prime.

\begin{lem}\label{LemmaB}
For $\gamma = (d c_1 , d c_2, \dots d c_m)$, we have
\begin{equation*}
G^{\gamma}_{\el}(q_1, \dots, q_m)
= G^{(c_1, \dots, c_m)}_{\floor{\el /d}}(q_1^{d}, \dots, q_m^d).
\end{equation*}
\end{lem}

\begin{proof}
It is straightforward.
\end{proof}

Let us calculate several examples of
$G^{\gamma}(q_1, \dots, q_m)$
for special parameters $\gamma$
which are rather non-trivial.

\begin{example}
Let $\gamma= (c,1)$ and calculate $G^{(c,1)}_{\el}(q_1,q_2)$.
Putting  $d=\floor{\el/ c}$, we see by Lemma \ref{LemmaA} that 
\begin{equation}\label{2variablesCase}
\begin{split}
G^{(c,1)}_{\el}(q_1,q_2)
& = \sum_{n=0}^d q_1^{cn} G^{1}_{c n} (q_2) \\
& = \sum_{n=0}^d q_1^{cn} \frac{1-q_2^{cn+1}}{1-q_2} 
= \frac1{1-q_2} \left( \frac{1-q_1^{c(d+1)} }{ 1-q_1^c}
- q_2 \frac{1-(q_1 q_2)^{c(d+1)}}{1-(q_1 q_2)^{c}} \right) \\
&= \frac{(1-q_1^{c(d+1)})(1-(q_1 q_2)^c)-q_2(1-q_1^c)(1-(q_1 q_2)^{c(d+1)})}
{(1-q_2)(1-q_1^c)(1-(q_1 q_2)^{c})}. 
\end{split}
\end{equation}
If we take a limit $\el \to \infty$,
we obtain
\begin{equation*}
G^{(c,1)}_{\infty}(q_1,q_2)
= \frac{1-q_2 +q_1^c q_2 -q_1^c q_2^c}{(1-q_2)(1-q_1^c)(1-q_1^c q_2^c)}.
\end{equation*}
Therefore we have 
\begin{equation*}
\begin{split}
Z^{(c,1)}_{\infty} (t_1,t_2)
&= \sum_{n_2|n_1^c} n_1^{-ct_1} n_2^{-t_2}\\
&= \zeta(t_2) \zeta(c t_1) \zeta(c(t_1+t_2))
\prod_{\prims} (1-p^{-t_2}+p^{-ct_1-t_2}-p^{-ct_1-ct_2})
\end{split}
\end{equation*}
for $\Re t_j >1 \ (j=1,2)$.
Thus we may have various possibility of Euler products of the form 
$\prod_{\prims} 
H(p^{-s},p^{-t})$, where $H(S,T)\in 1+S\mathbb \cdot C[S,T]+T\mathbb \cdot C[S,T]$ 
is arising from $Z^{(c,1)}_{\infty} (t_1,t_2)$. 
\qed
\end{example}

\begin{example}
We calculate $G^{(cd,c,1)}_{\infty} (q,q,q)$.
If we set $q_1=q_2=q$ in \eqref{2variablesCase},
then we have
\begin{equation*}
G^{(c,1)}_{\el}(q,q)
= \frac{(1-q^{c(d+1)})(1-q+q^c-q^{cd+c+1})}{(1-q)(1-q^{2c})}.
\end{equation*}
It follows from Lemma \ref{LemmaA} that
\begin{align*}
G^{(cd,c,1)}_{\infty} (q,q,q) 
&= \sum_{n=0}^{\infty} q^{cdn} G^{(c,1)}_{cdn}(q,q)
= \sum_{n=0}^{\infty} q^{cdn}
\frac{(1-q^{c(dn+1)})(1-q+q^c-q^{cdn+c+1})}{(1-q)(1-q^{2c})} \\
&= \frac{1}{(1-q)(1-q^{2c})}
\left\{ \frac{1-q+q^c}{1- q^{cd}} 
-\frac{q^c(1+q^c)}{1-q^{2cd}} + \frac{q^{2c+1}}{1- q^{3cd}}
\right\}.
\end{align*}
In particular, if we put $d=1$ then 
\begin{equation*}
G^{(c,c,1)}_{\infty}(q,q,q) 
= \frac{1-q+q^c-q^{c+1}+q^{2c+1}-q^{3c}+q^{3c+1} -q^{4c}}
{(1-q)(1-q^{2c})^2(1-q^{3c})} 
= \frac{1-q+q^c-q^{c+1}+q^{2c}}{(1-q)(1-q^{2c})(1-q^{3c})}.
\end{equation*}
Hence we obtain
\begin{equation*}
Z^{(c,c,1)}_{\infty}(s)
= \sum_{n_3|n_2^c|n_1^c} (n_1^c n_2^c n_3)^{-s}
= \zeta(s) \zeta(2cs) \zeta(3cs)
\prod_{\prims} (1-p^{-s} +p^{-cs} - p^{-(c+1)s} + p^{-2cs}).
\end{equation*}
\qed
\end{example}

\section{Multiple zeta functions and powerful numbers}

In this section, we study $Z^{\gamma}_\infty (s)$ for $\gamma=(k,k,\dots,k,1)$
in connection with a certain generalized notion of  \term{powerful numbers}.

Let us first recall the definition of powerful numbers. A positive number
$n$ is called a $k$-powerful number if 
$\ord_p n \geq k$ for any prime number $p$ unless $\ord_p n=0$
(see, e.g. \cite{IvicShiu}, \cite{Ivic}). Extending this we arrive at a
new notion,  
\term{$\el$-step $k$-powerful numbers}; 
a positive integer $n$ is said to be an $\el$-step $k$-powerful number 
if $n$ satisfies the condition that 
$\ord_p n= 0,k,2k,\dots,(\el-1)k$
or $\ord_p n \geq \el k$ for any prime number $p$. 
Clearly, if $n$ is an $\el$-step $k$-powerful number,
then $n$ is again a $j$-step $k$-powerful number
for each $j\;(1 \leq j \leq \el)$. 
In particular, $1$-step $k$-powerful numbers are nothing but
the usual $k$-powerful numbers. Note also that every natural number 
is an $\el$-step 1-powerful number for any $\el$; this agrees with the 
claim for $k=1$ in Theorem \ref{Fkl} below. 

As an example of $\el$-step $k$-powerful numbers, 
for instance, the first few of 2-step 2-powerful numbers are 
$1, 4, 9, 16, 25, 
32, 36, 49, 64, 81, 100, 121, 128, 144, 169,196, 225, 243, \ldots$. 
We note that in general, an $\el$-step $k$-powerful number $n$ has the
following canonical representation:
$$
n=a_1^ka_2^{2k}\cdots a_{\el}^{(\el-1)k}\times m,
$$
where $a_1,\cdots, a_\el$ are square-free, $m$ is  $\el k$-powerful
and these satisfy $\gcd(a_1,\cdots, a_\el, m)=1$. 
Note that here a $k$-powerful number $m$ is uniquely
expressed as $m=b_1^kb_2^{k+1}\cdots b_k^{2k-1}$ if we stipulate that
$b_2,\cdots b_k$ are all square-free. 

Let us put
\begin{equation*}
f_{k,\el}(n)\deq
\begin{cases}
1 & \text{if $n$ is an $\el$-step $k$-powerful number}, \\
0 & \text{otherwise}
\end{cases}
\end{equation*}
for a positive integer $n$. We define also 
$F_{k,\el}(s) \deq \sum_{n=1}^{\infty} f_{k,\el}(n) n^{-s}$.
This arithmetic function $f_{k,\el}(n)$ is multiplicative
with respect to $n$.
Note that $F_{1,\el}(s)=\zeta(s)$ for any $\el$.
We show that $Z^\gamma_\infty(s)$ is represented by the product
of the Riemann zeta functions times $F_{k,\el}(s)$.
Namely, we have the

\begin{thm}\label{ThmPowerful}
Let $k, \el$ be positive integers
and put $\gamma=(\overbrace{k,k,\dots,k,}^\el 1)$.
Then we have
\begin{equation}
Z^{\gamma}_\infty(s)
= \sum_{n_{\el+1}| n_{\el}^k | \cdots | n_1^k}
(n_1^k \cdots n_{\el}^k n_{\el+1})^{-s}
= F_{k,\el}(s) \prod_{j=2}^{\el+1} \zeta(jks) \qquad (\Re s>1).
\end{equation}
\end{thm}

\begin{remark}
When $k=1$, this theorem gives Corollary \ref{MMzetaInfty}.
\qed
\end{remark}

To prove the theorem,
we prepare the following two lemmas.

\begin{lem}\label{GeneratingPowerful}
Let $\gamma=(\overbrace{k,k, \dots, k,}^\el 1)$.
Then we have
\begin{equation*}
G^{\gamma}_\infty (q)
= G^{\gamma}_{\infty} (\overbrace{q,q, \dots,q}^{\el+1})
= \prod_{j=1}^{\el+1}
\frac{1}{1-q^{jk}} \cdot \frac{1-q+q^{\el k+1}-q^{k(\el+1)}}{1-q}.
\end{equation*}
\end{lem}

\begin{proof}
By definition, it follows that
\begin{align*}
 G^{(k,\dots,k,1)}_\infty (q)
&= \sum_{\begin{subarray}{c}
\lambda_1 \ge \dots \ge \lambda_{\el} \ge \lambda_{\el+1} \ge 0 \\
k|\lambda_j \ (1 \leq j \le \el)
\end{subarray}}
q^{\lambda_1+ \cdots + \lambda_{\el} + \lambda_{\el +1}}
= \sum_{n=0}^\infty q^n
\sum_{\begin{subarray}{c}
\lambda_1 \ge \dots \ge \lambda_{\el} \ge n \\
k|\lambda_j \ (1 \leq j \leq \el)
\end{subarray}}
q^{\lambda_1 + \cdots + \lambda_\el} \\
&= \sum_{\begin{subarray}{c}
\lambda_1 \geq \dots \geq \lambda_{\el} \geq 0 \\
k|\lambda_j \ (1 \leq j \leq \el)
\end{subarray}}
q^{\lambda_1 + \cdots + \lambda_\el}
+ \sum_{a=0}^\infty \sum_{b=1}^k q^{ak+b} \sum_{
\begin{subarray}{c}
\lambda_1 \geq \dots \geq \lambda_{\el} \geq ak+b \\ 
k|\lambda_j \ (1 \leq j \leq \el)
\end{subarray}}
q^{\lambda_1 + \cdots + \lambda_\el} \\
&= \sum_{\mu_1 \ge \dots \ge \mu_\el \ge 0}
q^{k( \mu_1 + \cdots + \mu_\el)}
+ \sum_{a=0}^\infty \sum_{b=1}^k q^{ak+b}
\sum_{\mu_1 \ge \cdots \ge \mu_\el \ge a+1}
q^{k (\mu_1+ \cdots +\mu_\el)} \\
&= G^\el_\infty(q^k) + \sum_{a=0}^\infty q^{ak}
\sum_{b=1}^k  q^{b} \sum_{\nu_1 \geq \cdots \geq \nu_\el \geq 0}
q^{k\{ (\nu_1+a+1)+ \cdots + (\nu_\el+a+1)\}}  \\
&= G^\el_\infty (q^k) + q^{\el k}
\left( \sum_{a=0}^\infty q^{k(\el+1)a} \right)
\left( \sum_{b=1}^k q^b \right) 
\left( \sum_{\nu_1 \geq \cdots \geq \nu_\el \geq 0}
q^{k(\nu_1+ \cdots +\nu_\el)} \right).
\end{align*}
Since $G^\el_\infty(q) = \prod_{j=1}^{\el} (1-q^j)^{-1}$,
we have
\begin{align*}
G^{(k,\dots,k,1)}_\infty (q) 
&= \prod_{j=1}^\el \frac{1}{1-q^{jk}}
\left( 1+ \frac{q^{k \el}}{1-q^{k(\el+1)}}
\frac{q(1-q^k)}{1-q} \right) \\
&= \prod_{j=1}^{\el+1} \frac1{1-q^{jk}}
\frac{(1-q)(1-q^{k(\el+1)}) + q^{\el k+1}(1-q^k)}{1-q} \\
&= \prod_{j=1}^{\el+1} \frac1{1-q^{jk}}
\frac{1-q+q^{\el k+1}-q^{k(\el+1)}}{1-q}.
\end{align*}
This proves the assertion. 
\end{proof}

\begin{lem}\label{ProductPowerful}
We have
\begin{equation*}
\frac{1-q+q^{\el k+1}-q^{k(\el+1)}}{1-q}
= (1-q^k) \left(1+q^k +q^{2k} + \cdots + q^{(\el-1)k}
+ q^{\el k} \sum_{j=0}^\infty q^j \right).
\end{equation*}
\end{lem}

\begin{proof}
The calculation is straightforward.
Actually, we have
\begin{align*}
& \quad (1-q^k)
\left( 1+q^k +q^{2k} + \cdots + q^{(\el-1)k} + q^{\el k}
\sum_{j=0}^\infty q^j \right)\\
& = (1-q^k)
\left( \frac{1-q^{\el k}}{1-q^k}
+ \frac{q^{\el k}}{1-q} \right) \\
& = \frac1{1-q}
\left\{ (1-q^{\el k})(1-q) + q^{\el k} (1-q^k) \right\}
= \frac{1-q+q^{\el k+1}-q^{k(\el+1)}}{1-q}.
\end{align*}
\end{proof}

\begin{proof}[Proof of Theorem \ref{ThmPowerful}]
It follows from Lemma \ref{GeneratingPowerful}
and Lemma \ref{ProductPowerful} that
\begin{align*}
Z^{(k,\dots, k,1)}_\infty (s) 
&= \prod_{\prims} G^{(k,\dots,k,1)}_\infty (p^{-s}) \\
&= \prod_{\prims}
\left( \prod_{j=2}^\infty \frac{1}{1-p^{-jks}} \right)\\
& \quad \times \prod_{\prims}
(1+ p^{-ks} + p^{-2ks} + \cdots + p^{-(\el-1)ks}
+ p^{-\el k s}+ p^{-(\el k+1)s} + \cdots ) \\
&= \prod_{j=2}^{\el+1} \zeta(jks) \cdot F_{k,\el}(s).
\end{align*}
This completes the proof of the theorem.
\end{proof}

Now we determine the condition if 
the Dirichlet series $F_{k,\el}(s)$ can be 
meromorphically extended to the whole complex plane $\bC$.
We recall the following Estermann's result \cite{Estermann}
(see \cite{K} for a generalization). A
polynomial
$f(T) \in 1+ T \cdot \bC [T]$ is said to be \term{unitary}
if and only if there is a unitary matrix $M$
such that $f(T) = \det (1-MT)$. 

\begin{lem}\label{Estermann}
For a polynomial $f(T) \in 1+ T \cdot \bC[T]$,
put $L(s,f) = \prod_{\prims} f(p^{-s})$.
Then
\begin{enumerate}
\item $f(T)$ is unitary
if and only if $L(s,f)$ can be extended
as a meromorphic function on $\bC$.
\item $f(T)$ is not unitary
if and only if $L(s,f)$ can be extended
as a meromorphic function in $\Re s >0$ with
the natural boundary $\Re s=0$ {\upshape ;}
each point on $\Re s=0$ is a limit-point of poles
of $L(s,f)$ in $\Re s>0$.
\end{enumerate}
\qed
\end{lem}

Since 
\begin{align*}
F_{k,\el}(s) 
&= \prod_{\prims} 
\left( 1+p^{-ks} + p^{-2ks} + \cdots + p^{-(\el-1)ks}
+ p^{-\el ks} \sum_{j=0}^{\infty} p^{-js} \right) \\
&= \zeta (s) \zeta(ks) \prod_{\prims}
(1-p^{-s}+p^{-(\el k+1)s} -p^{-k(\el+1)s})
\end{align*}
by Lemma \ref{ProductPowerful},
we have only to see whether the polynomial
$G_{k,\el}(T) \deq 1-T+T^{\el k+1} -T^{k(\el+1)}$
is unitary or not.
The polynomial $G_{k,\el}(T)$ can be expressed as
$G_{k,\el}(T) = (1-T^k) H_{k,\el}(T)$
with
$H_{k,\el}(T) \deq 1+(T^k-T)\sum_{j=0}^{\el-1} T^{kj}$.

\begin{prop}\label{Unitarity}
Let $k$ and $\el$ be positive integers.
The polynomial $G_{k,\el}(T)$ is unitary if and only if $k=1, 2$.
\end{prop}

In order to prove this proposition,
we need the following two lemmas.

\begin{lem}\label{RootG(T)}
Let $k$ be a positive integer such that $k \geq 3$.
Then the unitary root $\alpha$ (i.e. $\abs{\alpha}=1$) 
of the polynomial $G_{k,\el}(T)$
must satisfy $\alpha^k=1$ or $\alpha^{k-2}=1$.
\end{lem}

\begin{proof}
Let $\alpha=e^{2 \pi i \theta} \not=1 \ (\theta \in \bR)$
be a unitary root of $G_{k,\el}(T)$.
Since $G_{k,\el}(T)/(1-T) = 1+T^{\el k+1} (1-T^{k-1})/(1-T)$,
we have $1+ \alpha^{\el k+1} (1-\alpha^{k-1})/(1-\alpha)=0$ so that
$\abs{(1-\alpha^{k-1})/(1-\alpha)}
= \abs{\alpha^{-(\el k+1)}}=1$.
Hence we see that 
$\Re \alpha^{k-1} = \Re \alpha$, that is,
$\cos 2\pi (k-1)\theta- \cos 2\pi \theta
=-2\sin \pi k \theta \sin \pi (k-2) \theta =0$. 
Thus we conclude that either $k \theta \in \bZ$ or $(k-2)\theta \in \bZ$.
This proves the lemma.
\end{proof}

\begin{lem}\label{RootMultiplicity}
Let $k$ be a positive integer such that $k \geq 3$. Suppose that 
a complex number $\alpha$ satisfies $\alpha^{k-2}=1$.
Then  $G''_{k,\el}(\alpha) \not=0$.
\end{lem}

\begin{proof}
Since
$G''_{k,\el}(T) = (\el k+1)\el k T^{\el(k-2) +2 \el-1}
- (k \el+k) (k \el+k-1) T^{(k-2)(\el+1)+2\el}$,
if we assume that $\alpha$ satisfies $G''_{k,\el}(\alpha)=0$
and $\alpha^{k-2}=1$, we have
$G''_{k,\el}(\alpha)= (\el k+1)\el k \alpha^{2 \el-1}
- (k \el+k) (k \el+k-1) \alpha^{2\el}=0$.
This immediately shows that 
$\alpha=\frac{\el (\el k+1)}{(\el+1)(k \el +k-1)}$,
which contradicts $\alpha^{k-2}=1$.
Hence the assertion follows.
\end{proof}

\begin{proof}[Proof of Proposition \ref{Unitarity}]
Let $\el$ be a positive integer.
Since the unitarity of $G_{k,\el}(T)$
and that of $H_{k,\el}(T)$ are equivalent, it suffices to 
check the unitarity of $H_{k,\el}(T)$.
It is clear that $H_{k, \el}(T)$ is a unitary polynomial
when $k=1, 2$.
Actually we have
$H_{1,\el}(T)=1$ and
$H_{2,\el}(T) = 1-T+T^{2\el+1}-T^{2\el+2} = (1-T)(1+T^{2\el+1})$,
which are indeed unitary.

Suppose $k \ge 3$.
If $G_{k,\el}(T)$ is unitary,
then every root of $H_{k,\el}(T)$ satisfies
$\alpha^{k}=1$ or $\alpha^{k-2}=1$
by Lemma \ref{RootG(T)}.
However, if $\alpha^k=1$ we immediately see that
$H_{k,\el}(\alpha)= 1+(1-\alpha)\el$. Thus, $H_{k,\el}(\alpha)$ can not be 
0 because of the unitarity of $\alpha$.
Thus any root of $H_{k,\el}(T)$ must satisfy
$\alpha^{k-2}=1$ and $\alpha^k \neq 1$.
By Lemma \ref{RootMultiplicity},
the multiplicity of these roots of $H_{k,\el}(T)$ is at most $2$.
Since $H_{k,\el}(T)$ is assumed to be unitary
and $H_{k,\el}(1) \not=0$,
it follows that $2(k-3) \geq \deg H_{k,\el}(T) = \el k$.
This is possible only when $\el = 1$.

Therefore it is enough to prove that
$H_{k,1}(T)$ is not unitary for $k \ge 3$.
We put $H_k(T)=H_{k,1}(T)=1-T+T^k$ for simplicity.

Suppose that $k$ is an odd integer such that $k \geq 3$. 
Then $H_k(T)$ has a real root in the interval $(-2,-1)$
since $H_k(-1)=1>0$ and $H_k(-2) = 3-2^k <0$.
This implies the polynomial $H_k(T)$ is not unitary.

Thus, only we have to consider is the case
where $k$ is even and  $k \geq 4$.
Suppose that $H_k(T)$ is unitary
and let $e^{i\theta}$ \ $(-\pi < \theta \leq \pi)$
be its unitary root. 
Then we see that $\theta$  satisfies the equations 
$\cos k\theta = \cos \theta -1$ and $\sin k \theta = \sin \theta$.
Since $1= \sin^2 k\theta+ \cos^2 k\theta = 2-2 \cos \theta$,
we have $\theta = \pm\pi/3$.
Further,
since $1= \cos^2 \theta + \sin^2 \theta
= (\cos k \theta +1)^2+ \sin^2 k\theta = 2 \cos k \theta+2$,
we have $\cos(k \pi/3) = -1/2$.
Hence we see that either $k \equiv 2 \; \text{or} \; 4 \pmod{6}$. 
On the other hand, since
$1= (\cos \theta - \cos k \theta)^2 + (\sin \theta - \sin k \theta)^2
= 2 - 2\cos (k-1)\theta$,
we have $\cos((k-1)\pi/3) = 1/2$.
It follows that either $k \equiv 0 \;\text{or} \; 2 \pmod{6}$ holds.
Thus we have $k \equiv 2 \pmod{6}$.

Now we show that every unitary root of $H_k(T)$ is simple.
If we assume that $\beta$ is a multiple root of $H_k(T)$,
it follows that
$\beta^k-\beta+1=0$ and $k \beta^{k-1}-1=0$. 
Then we have $\abs{\beta}= k^{-1/(k-1)}$
and $\beta=k \beta^k$ by the second equation.
On the other hand,
by the first equation,
we obtain $1= \beta-\beta^k = (k-1)\beta^k$
so that $\abs{\beta}=(k-1)^{-1/k}$.
Therefore we have $k^k= \abs{\beta}^{-k(k-1)} = (k-1)^{k-1}$,
but this contradicts to the assumption of the unitarity of $H_k(T)$.
This completes the proof of the proposition.
\end{proof}

Finally we obtain the following theorem, which gives 
actually a generalization of the result in \cite{IvicShiu} concerning the
powerful numbers. 

\begin{thm}\label{Fkl}
Let $k$ and $\el$ be positive integers. 
When $k=1, 2$ we have 
\begin{equation*}
F_{1,\el}(s)=\zeta(s),\qquad
F_{2,\el}(s)= \frac{\zeta(2s) \zeta((2\el+1)s)}{\zeta(2(2\el+1)s)}. 
\end{equation*}
When $k \geq 3$, $F_{k,\el}(s)$ can be extended
as a meromorphic function in $\Re s>0$ and has a natural boundary $\Re s=0$.
\end{thm}

\begin{proof}
The assertion is clear from Lemma \ref{Estermann}
and Proposition \ref{Unitarity}.
\end{proof}
As a consequence we have the following result. 
\begin{cor}
Let $Z^{(k,\dots,k,1)}_{\infty}(s)$ be
as in Theorem \ref{ThmPowerful}.
Then we have for $k=1, 2$ 
\begin{equation*}
Z^{(1,\dots,1,1)}_{\infty}(s)= \prod_{j=1}^{\el+1} \zeta(js),\qquad
Z^{(2,\dots,2,1)}_\infty (s)
=\frac{\zeta((2\el+1)s)}{\zeta(2(2\el+1)s)}
\prod_{j=1}^{\el+1} \zeta(2js),
\end{equation*}
and for $k \geq 3$ the function 
$Z^{(k,\dots,k,1)}_{\infty}(s)$ can be meromorphically extended
to the half plane $\Re s>0$ 
with a natural boundary $\Re s=0$.
\qed
\end{cor}

\section{Closing remarks}\label{ConcludingRemarks}

We make two small remarks;
the first one is to describe the relation
between a multiple finite Riemann zeta function
and the theory of elementary divisors,
and the second one is to examine a multiple Eisenstein series defined
via the multiple finite Riemann zeta function.

\medskip

\noindent
$\bullet$
The isomorphism classes of abelian groups $A$ of order $n$
are parametrized by the map $\blambda$
from the set of all prime numbers to that of partitions
such that $n= \prod_{\prims} p^{\abs{\blambda(p)}}$ and
$$
A \cong \bigoplus_{\prims} \bigoplus_{j=1}^{\ell(\blambda(p))}
\bZ / p^{\blambda_j(p)}\bZ,
$$
where $\abs{\blambda(p)}$ and $\length{\blambda(p)}$ are
the size and the length of the partition
$\blambda(p)= (\blambda_j(p))_{j \ge 1}$ respectively.
The multiple finite Riemann zeta function $Z^m_N(s)$ is expressed also as
$Z^m_N(s) = \sum_{n| N^m} g^m_N(n) n^{-s}$.
Here $g^m_N(n)$ denotes the number of isomorphism classes
of abelian groups of order $n$,
parametrized by $\blambda$ such that
$\blambda_1(p) \leq m$ and $\ell(\blambda(p)) \leq \ord_p N$
for all $p$.
It is clear that $g^m_N(n)$ is multiplicative
with respect to $n$ and $N$.
If we put
$g^m_{\infty}(n) \deq \lim_{N \to \infty} g^m_{N}(n)$,
then $g^m_{\infty}(n)$ equals the number of the isomorphism classes
of abelian groups $A$ of order $n$
which is the direct sum of $p$-groups $A_p$
such that $p^m A_p =0$ for $p|n$.

We now just mention about the asymptotic average for $g^m_\infty(n)$
and $Z^m_n(\sigma)$ ($\sigma \in \bR$)
with respect to $n$.
Thus we need the Tauberian theorem below (see, e.g. \cite{MM}).

\begin{lem}\label{Tauberian}
Let $F(t)=\sum^{\infty}_{n=1}a_n n^{-t}$ be
a Dirichlet series with non-negative real coefficients
which converges absolutely for $\Re(t)>\beta$.
Suppose that $F(t)$ has a meromorphic continuation
to the region $\Re(t)\ge \beta$ with a pole
of order $\alpha+1$ at $t=\beta$ for some $\alpha \ge 0$.
Put
\begin{equation*}
c\deq\frac{1}{\alpha!}\,
\lim_{t \to \beta}(t-\beta)^{\alpha+1}F(t).
\end{equation*}
Then we have
\begin{equation*}
 \sum_{n\le x}a_n=(c+o(1))x^{\beta} (\log x)^{\alpha}
\end{equation*}
as $x \to \infty$.
\qed
\end{lem}

Using the lemma above, we have easily the following facts: 
Let $m$ be a positive integer.
\begin{enumerate} 
\item
\begin{equation}\label{asy1}
\sum_{n \leq x} g^m_{\infty}(n)
= (\zeta(2) \zeta(3) \dotsb \zeta(m)+o(1)) x
\end{equation}
as $x \to \infty$.
In other words, the asymptotic average of $g^m_\infty (n)$
with respect to $n$ is given by
$\zeta(2) \zeta(3) \dotsb \zeta(m)$.
\item
For a fixed $\sigma>0$,
\begin{align}
\sum_{n \leq x} Z^m_n(\sigma)
&= (\zeta(\sigma+1) \zeta(2\sigma+1) \cdots \zeta(m\sigma+1)+o(1)) x, \\
\sum_{n \leq x} Z^m_n(-\sigma)
&= (\zeta(\sigma+1) \zeta(2\sigma+1) \cdots \zeta(m\sigma+1) +o(1)) x^{1+m \sigma}, \\
\sum_{n \leq x} Z^m_n(0) 
&=\sum_{n \leq x} \prod_{\prims} \binom{\ord_p n+m}{m}
= \left( \frac{1}{m!} +o(1) \right)x (\log x)^m \label{asy}
\end{align}
as $x \to \infty$.
\end{enumerate}
Actually,
since $\zeta(s)$ has a single pole at $s=1$
and $\Res_{s=1} \zeta(s) = 1$,
it follows the first formula \eqref{asy1} from  \eqref{m-pleZeta}
and Lemma \ref{Tauberian}.
Next, fix $\sigma \in \bR$.
By \eqref{m-pleZeta(s,t)},
we have
$\zeta^m (\sigma, t)
= \sum_{n=1}^{\infty} Z^m_n (\sigma) n^{-t}
= \prod_{k=0}^m \zeta(t+ k\sigma)$.
This shows that the abscissa of absolute convergence of the  
Dirichlet series $\zeta^m (\sigma, t)$  
is given by 
$t =\max\{1, 1-m \sigma \}$.
Hence the remaining formulas follow similarly.

\medskip

\begin{remark}
Put $g(n)\deq\lim_{m \to \infty} g^m_\infty (n)$.
It is well-known that
$\sum_{n \leq x} g(n)
=\Big( \prod_{k=2}^{\infty}\zeta(k)\Big) x + O(\sqrt{x})$.
\qed
\end{remark}

\medskip

\begin{remark}
Since $Z^1_n(0) = d(n)\deq\sum_{d|n}1$,
it follows from \eqref{asy}
that $\sum_{n \leq x} d(n) \sim x \log x$.
It is also well-known \cite{Z} that
there exists a constant $C$ such that
$\sum_{n \leq x} d(n) = x \log x + C x +O(\sqrt{x})$
in an elementary way.
See, e.g. \cite{A}.
\qed
\end{remark}

\medskip

\noindent
$\bullet$ We define a multiple Eisenstein series
with parameter $s$ of type $m$ by
\begin{equation*}
E_s^m(q)=\sum_{n=1}^\infty Z_n^m(1-s)q^{n}.
\end{equation*}
We sometimes write $E_s^m(\tau)$ instead of $E_s^m(q)$
when $q=e^{2\pi i \tau}$ with $\tau\in\bC,\;\Im(\tau)>0$.
It is obvious that
$E_k^1(q)$ is (essentially) the usual holomorphic
Eisenstein series of weight $k$.
In this remark we make an experimental study
of this multiple Eisenstein series.
Here we assume that $m=2$;
we shall make a detailed study of these series in somewhere else.

First we observe the following simple relation.

\begin{lem}\label{Eisen1} 
We have
\begin{equation*}
E_{s+1}^2(q)
=\sum_{\el=1}^\infty \sum_{N=1}^\infty \sigma_s(N)N^s q^{N\el}.
\end{equation*}
\end{lem}

\begin{proof}
The calculation is straightforward. Actually, one has
\begin{align*}
& \sum_{\el=1}^\infty
\sum_{N=1}^\infty \sigma_s(N)N^s q^{N\el}
=\sum_{\el=1}^\infty
\sum_{m=1}^\infty \sum_{n=1}^\infty (nm)^sn^s q^{nm\el}
=\sum_{\el=1}^\infty
\sum_{j=1}^\infty \sum_{n|j} j^s n^s q^{j\el}\\
= &\sum_{N=1}^\infty
\sum_{n|j}\sum_{j|N}n^s j^s q^N
= \sum_{N=1}^\infty Z_N^2(-s)q^N
=  E_{s+1}^2(q).
\end{align*}
This proves the assertion.
\end{proof}

Recall now the Fourier expansion of
the holomorphic Eisenstein series $E_{k+1}(\tau)$
of weight $k+1$ with $k$ being odd;
\begin{align*}
E_{k+1}(\tau)
&=\frac12\sum_{(c,d)=1}(c\tau+d)^{-k-1}
=1+\frac1{\zeta(k+1)}\sum_{m=1}^\infty
\sum_{n\in \bZ}(m\tau+n)^{-k}\\
&=1+\frac1{\zeta(k+1)}\frac{(2\pi i)^{k+1}}{k!}
\sum_{n=1}^\infty \sigma_{k}(n)q^n\;
\left(=1+\frac1{\zeta(k+1)}
\frac{(2\pi i)^{k+1}}{k!}E^1_{k+1}(\tau)\right).
\end{align*} 
Taking $k$-times derivative of $E_{k+1}(\tau)$, we have easily
\begin{equation*}
\frac{d^{k}}{d\tau^{k}}E_{k+1}(\tau)=
\frac{(2\pi i)^{2k+1}}{\zeta(k+1) k!}
\sum_{n=1}^\infty \sigma_{k}(n)n^kq^n.
\end{equation*}
Hence by Lemma \ref{Eisen1},
we immediately obtain the expression of $E^2_{k+1} (\tau)$ when $k$ is odd:
\begin{equation}
E_{k+1}^2(\tau) = \frac{\zeta(k+1) k!}{(2\pi i)^{2k+1}}
\sum_{\el=1}^\infty
\left(\frac{d^{k}}{d\tau^{k}}E_{k+1}\right)(\el\tau).
\end{equation}

\medskip

\begin{remark}
There is another expression of $E_{k+1}^2(\tau)$ for $k$ being odd 
as follows.
\begin{equation}
E_{k+1}^2(\tau)=- \frac{(2k)!}{(2\pi i)^{2k+1}}
\sum_{\el=1}^\infty \sum_{\begin{subarray}{c}
(c,d)=1\\
c>0
\end{subarray}}
\sigma_k(c\el)\el^{-2k-1}(c\tau+d)^{-2k-1}.
\end{equation}
\qed
\end{remark}

\par\noindent
$\bullet$ It is clearly interesting to make an extensive 
 study of the multi-variable 
and multi-parameter version of multiple finite Riemann zeta functions from 
the symmetric functions point of view.  


\noindent
\textsc{\by{KIMOTO}{Kazufumi}}\\
Department of Mathematical Science, University of the Ryukyus.\\
Nishihara, Okinawa, 903-0213 JAPAN.\\
\texttt{kimoto@math.u-ryukyu.ac.jp}\\

\noindent
\textsc{\by{KUROKAWA}{Nobushige}}\\
Department of Mathematics, Tokyo Institute of Technology.\\
Oh-okayama Meguro-ku, Tokyo, 152-0033 JAPAN.\\
\texttt{kurokawa@math.titech.ac.jp}\\

\noindent
\textsc{\by{MATSUMOTO}{Sho}}\\
Graduate School of Mathematics, Kyushu University.\\
Hakozaki Higashi-ku, Fukuoka, 812-8581 JAPAN.\\
\texttt{ma203029@math.kyushu-u.ac.jp}\\

\noindent
\textsc{\by{WAKAYAMA}{Masato}}\\
Faculty of Mathematics, Kyushu University.\\
Hakozaki Higashi-ku, Fukuoka, 812-8581 JAPAN.\\
\texttt{wakayama@math.kyushu-u.ac.jp}

\end{document}